\documentclass[english,a4paper]{amsart}
\usepackage[T1]{fontenc}
\usepackage{ae}
\usepackage{aecompl}
\usepackage{amssymb}
\usepackage{mathrsfs}

\input{xy}
\xyoption{all}
\SelectTips{cm}{}
\usepackage{babel}
\usepackage[pdftex]{graphicx}
\usepackage[leqno]{amsmath}
\usepackage{amsthm}

\usepackage[svgnames]{xcolor}
\usepackage{tikz}
\usetikzlibrary{shapes,patterns,calc,snakes}

\renewcommand\phi{\varphi}
\renewcommand\epsilon{\varepsilon}
\renewcommand\theta{\vartheta}

\newcommand\mbb{\mathbb}

\newcommand\R{\mbb{R}}

 \newcommand{\la}[0]{\lambda}

\DeclareMathOperator\conv{conv}

\DeclareMathOperator\cc{cc}
\numberwithin{equation}{section}

\theoremstyle{plain}

\newtheorem{Prop}[equation]{Proposition}
\newtheorem{Cor}[equation]{Corollary}

\newtheorem*{Thm*}{Theorem}
\newtheorem*{Satz*}{Satz}
\newtheorem*{Prop*}{Proposition}
\newtheorem*{Cor*}{Corollary}
\newtheorem*{Lemma*}{Lemma}
\newtheorem*{Hilfssatz*}{Lemma}
\newtheorem*{Sublemma*}{Sublemma}
\newtheorem*{Conjecture*}{Conjecture}

\theoremstyle{definition}

\newtheorem{Ex}[equation]{Example}

\newtheorem{Rem}[equation]{Remark}

\newtheorem*{Def*}{Definition}
\newtheorem*{Defs*}{Definitions}
\newtheorem*{Ex*}{Example}
\newtheorem*{Exs*}{Examples}
\newtheorem*{LemmaDef*}{Lemma and Definition}
\newtheorem*{Notation*}{Notation}
\newtheorem*{Problem*}{Problem}
\newtheorem*{Question*}{Question}
\newtheorem*{Rem*}{Remark}
\newtheorem*{Rems*}{Remarks}
\newtheorem*{Warning*}{Warning}

\begin{document}
\title{A Note on the Convex Hull of finitely many Projections of Spectrahedra }
\author{Tim Netzer}
\address{Fachbereich Mathematik, Universit{\"a}t Konstanz, 78457 Konstanz, Germany}
\email{tim.netzer@uni-konstanz.de}

\author{Rainer Sinn}
\address{Fachbereich Mathematik, Universit{\"a}t Konstanz, 78457 Konstanz, Germany}
\email{rainer.sinn@uni-konstanz.de}

\date{\today}

\begin{abstract} A spectrahedron is a set defined by a linear matrix inequality. A projection of a spectrahedron is often  called a \textit{semidefinitely representable set}. We show that the convex hull of a finite union of such projections is again a projection of a spectrahedron. This improves upon the result of Helton and Nie \cite{HeltonNieNecSuffSDP}, who prove the same result in the case of bounded sets.\end{abstract}

\maketitle
\section{Introduction} Let $A_0,A_1,\ldots,A_n$ be real symmetric $k\times k$ matrices. The set $$\left\{ (x_1,\ldots, x_n)\in\R^n\mid A_0 +x_1A_1+ \cdots +x_nA_n\succeq 0\right\},$$ where $\succeq 0$ means positive semidefiniteness, is called a \textit{spectrahedron}. Spectrahedra are generalizations of polyhedra and occur as feasible sets for semidefinite optimization. 

A projection of a spectrahedron to a subspace of $\R^n$ is often called a \textit{semidefinitely representable set}.  Helton and Nie \cite{HeltonNieNecSuffSDP} conjecture that every convex semialgebraic set is such a projection.
See for example \cite{MR1284712, HeltonNieSDPrepr,HeltonNieNecSuffSDP,MR2292953,LasserreConvSets,NeSDP,NePlSch} for more detailed information on spectrahedra and their projections.

We prove that the convex hull of finitely many projections of spectrahedra is again a projection of a spectrahedron. This generalizes Theorem 2.2 from Helton and Nie \cite{HeltonNieNecSuffSDP}, which is the same result in the case that all sets are bounded or that the convex hull is closed. 
\section{Result}

\begin{Prop} If $S\subseteq \R^n$ is a projection of a spectrahedron, then so is $\cc(S)$, the conic hull of $S$.
\end{Prop}
\begin{proof} Since $S$ is a projection of a spectrahedron we can write $$S=\left\{x\in\R^n\mid \exists z\in\R^m\colon A+\sum_{i=1}^nx_iB_i +\sum_{j=1}^m z_jC_j\succeq 0\right\},$$ with suitable real symmetric $k\times k$-matrices $A,B_i,C_j$. Then with \begin{align*}C:=\{ x\in\R^n\mid & \exists \la, r \in\R, z\in\R^m\colon  \la A+\sum_{i=1}^n x_iB_i +\sum_{j=1}^m z_jC_j \succeq 0\  \wedge \\ & \quad \bigwedge_{i=1}^n \left(\begin{array}{cc}\la & x_i \\x_i & r\end{array}\right)\succeq 0\}\end{align*} we have $C=\cc(S)$ (note that $C$ is a projection of a spectrahedron, since the conjunction can be eliminated, using block matrices). 

To see "$\subseteq$" let some $x$ fulfill all the conditions from $C$, first with some $\la>0$. Then $a:=\frac{1}{\lambda}\cdot x$ belongs to $S$, using the first condition only. Since $x=\la\cdot a$, $x\in\cc(S)$. If $x$ fulfills the conditions with $\la=0$, then $x=0$, by the last $n$ conditions in the definition of $C$. So clearly also $x\in\cc(S)$.

For "$\supseteq$" take $x\in \cc(S)$. If $x\neq 0$ then there is some $\la>0$ and $a\in S$ with $x=\la a$. Now there is some $z\in\R^m$ with $A+\sum_i a_iB_i +\sum_j z_jC_j\succeq 0$. Multiplying this equation with $\la$ shows that $x$ fulfills the first condition in the definition of $C$. But since $\la> 0$, the other conditions can clearly also be satisfied with some big enough $r$. So $x$ belongs to $C$. Finally, $x=0$ belongs to $C,$ too.
\end{proof}

\begin{Rem}
The additional $n$ conditions in the definition of $C$ avoid problems that could occur in the case $\la=0.$ This is the main difference to the approach of Helton and Nie in \cite{HeltonNieNecSuffSDP}.
\end{Rem}

\begin{Cor} If $S_1,\ldots, S_t\subseteq\R^n$ are projections of spectrahedra, then also the convex hull $\conv(S_1\cup\cdots\cup S_t)$ is a projection of a spectrahedron.
\end{Cor}
\begin{proof} Consider $\widetilde{S}_i:=S_i\times\{1\}\subseteq\R^{n+1}$, and let $K_i$ denote the conic hull of $\widetilde{S}_i$ in $\R^{n+1}$. All $\widetilde{S}_i$ and therefore all $K_i$ are projections of spectrahedra, and thus the Minkowski sum $K:=K_1+\cdots + K_t$ is also  such a projection. Now  one easily checks $$\conv(S_1\cup\cdots\cup S_t)=\left\{x\in\R^n\mid (x,1)\in K\right\},$$ which proves the result.
\end{proof}

\begin{Ex} Let $S_1:=\{ (x,y)\in \R^2\mid x\geq 0, y\geq 0, xy\geq 1\}$ and $S_2=\{(0,0)\}$. Both subsets of $\R^2$ are spectrahedra, so the convex hull of their union, $$\conv(S_1\cup S_2)= \{ (x,y)\in\R^2\mid x>0, y>0 \} \cup \{ (0,0)\} ,$$ is a projection of a spectrahedron.

\end{Ex}

{\linespread{1}\bibliographystyle{dpbib} \bibliography{references}}

\end{document}